\documentclass[12pt]{amsart}
\usepackage{latexsym}

\calclayout

\renewcommand{\:}{\colon}
\newcommand{\+}{\nobreakdash-}
\renewcommand{\.}{\text{$\mskip .5\thinmuskip$}}
\renewcommand{\l}{\ell}
\newcommand{\ds}{\dots}
\newcommand{\dsb}{\dotsb}

\newcommand{\newcommandDOTSB}[2]
      {\newcommand{#1}{}\def#1{\DOTSB #2}}
\newcommand{\tinyskip}[1]{\mskip #1\thinmuskip}
\newcommand{\adjskip}{\tinyskip{.66}}
\newcommandDOTSB{\sub}{\subset}
\newcommandDOTSB{\ot}{\otimes}
\newcommandDOTSB{\op}{\oplus}
\newcommandDOTSB{\rarrow}{\longrightarrow}
\newcommandDOTSB{\larrow}{\longleftarrow}
\newcommandDOTSB{\lrarrow}
  {\adjskip\relbar\joinrel\relbar\joinrel\rightarrow\adjskip}
\newcommandDOTSB{\lorarrow}
  {\adjskip\relbar\joinrel\relbar\joinrel\relbar\joinrel
  \rightarrow\adjskip}
\newcommandDOTSB{\longrarrow}
  {\adjskip\relbar\joinrel\relbar\joinrel\relbar\joinrel
  \relbar\joinrel\rightarrow\adjskip}
\newcommandDOTSB{\lolongrarrow}
  {\adjskip\relbar\joinrel\relbar\joinrel\relbar\joinrel
  \relbar\joinrel\relbar\joinrel\rightarrow\adjskip}
\newcommandDOTSB{\longlongrarrow}
  {\adjskip\relbar\joinrel\relbar\joinrel\relbar\joinrel
  \relbar\joinrel\relbar\joinrel\relbar\joinrel
  \rightarrow\adjskip}
\newcommandDOTSB{\Longrarrow}
  {\adjskip\relbar\joinrel\relbar\joinrel\relbar\joinrel
  \relbar\joinrel\relbar\joinrel\relbar\joinrel
  \relbar\joinrel\relbar\joinrel\relbar\joinrel
  \relbar\joinrel\relbar\joinrel\rightarrow\adjskip}
\newcommandDOTSB{\Lolongrarrow}
  {\adjskip\relbar\joinrel\relbar\joinrel\relbar\joinrel
  \relbar\joinrel\relbar\joinrel\relbar\joinrel
  \relbar\joinrel\relbar\joinrel\relbar\joinrel
  \relbar\joinrel\relbar\joinrel\relbar\joinrel
  \rightarrow\adjskip}
\newcommandDOTSB{\Longlongrarrow}
  {\adjskip\relbar\joinrel\relbar\joinrel
  \relbar\joinrel\relbar\joinrel\relbar\joinrel
  \relbar\joinrel\relbar\joinrel\relbar\joinrel
  \relbar\joinrel\relbar\joinrel\relbar\joinrel
  \relbar\joinrel\relbar\joinrel\relbar\joinrel
  \relbar\joinrel\relbar\joinrel\rightarrow\adjskip}
\newcommandDOTSB{\llarrow}
  {\adjskip\leftarrow\joinrel\relbar\joinrel\relbar\adjskip}

\newcommand{\Z}{{\mathbb Z}}
\let\SSS\S
\renewcommand{\S}{\Sigma}
\renewcommand{\d}{\partial}

\DeclareMathOperator{\coker}{coker}
\DeclareMathOperator{\im}{im}
\DeclareMathOperator{\pr}{pr}
\DeclareMathOperator{\id}{id}
\DeclareMathOperator{\res}{res}
\DeclareMathOperator{\cor}{cor}
\DeclareMathOperator{\tr}{\mathit{tr}}
\DeclareMathOperator{\Hom}{Hom}
\DeclareMathOperator{\Ext}{Ext}
\DeclareMathOperator{\chr}{char}
\DeclareMathOperator{\Gal}{Gal}
\DeclareMathOperator{\Ann}{Ann}

\global\let\le\undefined
\global\let\ge\undefined
\DeclareMathSymbol{\le}{\mathrel}{AMSa}{"36}      
\DeclareMathSymbol{\ge}{\mathrel}{AMSa}{"3E}      

\newcommand{\oF}{\,\overline{\!F}}
\newcommand{\KM}{K^{\mathrm{M}}}

\newcommand{\Section}[1]{\bigskip\section{#1}\medskip}

\newtheoremstyle{sftheorem}
   {\topsep}{\topsep}{\slshape}{\parindent}{\sf}{.}{ }{}
\newtheoremstyle{sfremark}
   {\topsep}{\topsep}{}{\parindent}{\sf}{.}{ }{}
\theoremstyle{sftheorem}
\newtheorem{theorem}{Theorem}
\newtheorem{proposition}[theorem]{Proposition}
\newtheorem{corollary}[theorem]{Corollary}
\newtheorem{lemma}[theorem]{Lemma}
\newtheorem{conjecture}[theorem]{Conjecture}
\theoremstyle{sfremark}
\newtheorem*{remark}{Remark}

\begin{document}

\vspace*{-.5cm}

\title{Galois cohomology of certain field extensions\\
      and the divisible case of Milnor--Kato conjecture}
\author{Leonid Positselski}

\address{Independent University of Moscow}
\email{posic@mccme.ru}

\maketitle
\vspace{1mm}

 Let $F$ be a field and $m\ge2$ be an integer not divisible by
the characteristic of~$F$.
 Consider the absolute Galois group $G_F=\Gal(\oF/F)$, where
$\oF$~denotes the (separable) algebraic closure of~$F$.
 The famous \emph{Milnor--Kato conjecture}~\cite{Mil, Kato}
claims that the natural homomorphism of graded rings
  $$
   \KM_n(F)\ot\Z/m\lrarrow H^n(G_F,\.\mu_m^{\ot n})
  $$
from the Milnor K\+theory of the field~$F$ modulo~$m$ to
the Galois cohomology of~$F$ with cyclotomic coefficients
is an isomorphism (here, as usually, $\mu_m$ denotes the group
of $m$\+roots of unity in~$\oF$).
 One can see that it suffices to verify this conjecture
in the case when $m$~is a prime number.

 V.~Voevodsky~\cite{Voev} proved this conjecture for $m$ equal
to a power of~$2$.
 In his paper he also outlined a general approach to the Milnor-Kato
conjecture for arbitrary~$m$.
 The first step of his argument~\cite[Sections~5\+-6]{Voev} dealt
with a prime number~$\l$, a field $F$ having no finite extensions of
degree prime to~$\l$, and an integer $n\ge1$ such that
the group $\KM_{n+1}(F)$ is $\l$\+divisible.
 Assuming that the conjecture holds in degree less or equal to~$n$
for $m=\l$ and any field containing~$F$, Voevodsky was proving
that $H^{n+1}(G_F,\.\Z/\l)=0$.
 The main goal of this paper is to give a simplified elementary
version of Voevodsky's proof of this step.
 In particular, we do not need to consider fields transcendental
over~$F$ and we make no use of motivic cohomology at all.
 Our main result is formulated as follows.

\begin{theorem}
 \.Let $\l$ be a prime number, $n\ge1$ be an integer, and
$F$~be a field of characteristic not equal to~$\l$,
having no finite extensions of degree prime to~$\l$.
 Suppose that the homomorphism $\KM_n(L)/\l\KM_n(L)
\rarrow H^n(G_L,\.\Z/\l)$ is an isomorphism for any
finite extension~$L$ of the field~$F$.
 Furthermore, assume that $\KM_{n+1}(F)=\l\KM_{n+1}(F)$.
 Then one has $H^{n+1}(G_F,\.\Z/\l)=\nobreak 0$.
\end{theorem}

 We deduce the above theorem from two results, one related
to the Milnor K\+theory of fields and the other to
the cohomology of Galois groups.
 The first one is but a version of a lemma of Suslin~\cite{Sus},
which was also used in~\cite{Voev}.
 We only rewrote it for the Milnor K\+theory groups modulo~$\l$.
 The second statement is essentially the ``relative conjecture''
of Bloch and Kato~\cite{BK}, which was proven for $n=2$ by
Merkurjev and Suslin~\cite[\SSS15]{MS} already.
 It is obvious for $m=2$ (when it holds for any pro-finite group).
 For odd primes~$\l$ and composite numbers $m>2$ we show that
the desired cohomological sequence is exact provided that certain
Bockstein homomorphisms vanish.
 This is the key point of this paper.

\begin{proposition}
 \.Let $m\ge2$ be an integer, $G$ be a pro-finite group,
and $H\sub G$ be a normal subgroup such that $\S=G/H$
is a cyclic group of order~$m$.
 Let $T_{m^2}$ be a free module over $\Z/m^2$ endowed with
an action of~$G$, and let $T_m=m T_{m^2}=T_{m^2}/m T_{m^2}$
be the corresponding $G$\+module over~$\Z/m$.
 Assume that for a certain $n\ge0$ both Bockstein homomorphisms
 $$
  H^n(G,\.T_m)\lrarrow H^{n+1}(G,\.T_m)
  \quad \text{ and } \quad
  H^n(H,\.T_m)\lrarrow H^{n+1}(H,\.T_m)
 $$
vanish.  Then the following sequence of cohomology groups is exact
 \begin{multline*}
   H^n(G,\.T_m) \overset{\res}\lrarrow H^n(H,\.T_m)_\S
   \overset{\cor}\lrarrow H^n(G,\.T_m)    \\
   \overset{\!u\,\cup\.}\lrarrow H^{n+1}(G,\.T_m)
   \overset{\res}\lrarrow H^{n+1}(H,\.T_m)^\S
   \overset{\cor}\lrarrow H^{n+1}(G,\.T_m),
 \end{multline*}
where $u\in H^1(G,\.\Z/m)$ is the class corresponding
to the character $G\rarrow\S\simeq \Z/m$.
 (The letter~$\Sigma$ in the upper or lower index denotes
the invariants or the coinvariants with respect to a natural action
of the group~$\Sigma$.)
\end{proposition}

 In the second half of this paper we apply the same techniques
to obtain further exact sequences of Galois cohomology for
cyclic, biquadratic, and dihedral field extensions.
 In particular, our method proves the biquadratic exact sequences
conjectured by Merkurjev--Tignol~\cite{MT} and Kahn~\cite{Kahn}.
 In addition, we introduce an extended version of the classical
Bass--Tate lemma~\cite{BT} and deduce some corollaries about
generators and relations of annihilator ideals in
Galois cohomology rings.

\smallskip

 I would like to express my gratitude to Vladimir Voevodsky
for posing the problem and for very helpful conversations
and explanations.
 He also always urged me to write down this proof.
 I am grateful to Alexander Vishik for numerous very useful
discussions and to Roman Bezrukavnikov for one helpful remark.
 Most of this work was done during my stay at
the Institute for Advanced Study in Princeton and the rest
at the Institut des Hautes \'Etudes Scientifiques in Paris,
Max-Planck-Intitut f\"ur Mathematik in Bonn and
the Independent University of Moscow.
 I wish to thank all these four institutions.
 Finally, I want to thank the referees for their valuable
suggestions.

\Section{Suslin's Lemma for Milnor K-theory modulo~\protect{$\l$}}

 The goal of this section is to deduce Theorem~1 from Proposition~2.
 This may be the least elementary part of this paper, to the extent
that the proof of Corollary~4 is based on the existence of transfer
homomorphisms on the Milnor K\+theory groups.
 We will only use transfer maps for extensions of degree~$\l$
of fields having no extensions of degree prime to~$\l$.
 The existence of transfers for such field extensions
was established in~\cite{BT}; their basic properties are
the base change and the projection formula.

 Let us fix a prime number~$\l$, a number $n\ge 1$, and a field~$F$
with $\chr F\ne\l$, having no finite extensions of degree prime to~$\l$.
 Introduce the notation $k_m(K)=\KM_m(K)/\l$ for the Milnor
K\+theory groups modulo~$\l$ of a field~$K$.
 For a field extension $K\sub L$ we will denote by
$$
 i_{L/K}\:k_m(K)\lrarrow k_m(L)
   \quad\text{and}\quad
 \tr_{L/K}\:k_m(L)\lrarrow k_m(K)
$$
the inclusion and transfer homomorphisms.

\begin{lemma}
 \.Assume that for any finite field extensions $F\sub K\sub L$
with $[L:K]=\l$ the sequence of Milnor K\+theory groups
 $$
  k_n(K)\overset{i}\lrarrow k_n(L)_\S
  \overset{\tr}\lrarrow k_n(K), \qquad \S=\Gal(L/K)
 $$
is exact.
 Further assume that $[E:F]=\l$ and the map
$\tr_{E/F}\:k_n(E)\rarrow k_n(F)$ is surjective.
 Then the sequence
 $$
  k_{n+1}(F)\overset{i}\lrarrow k_{n+1}(E)_\S
  \overset{\tr}\lrarrow k_{n+1}(F), \qquad \S=\Gal(E/F)
 $$
is also exact.
\end{lemma}

\begin{proof}
 In order to prove that the map
$\tr_{E/F}\:k_{n+1}(E)_\S/i_{E/F}k_{n+1}(F)\rarrow k_{n+1}(F)$
is injective, we will construct a surjective homomorphism
 $$
  \phi\:k_{n+1}(F)\lrarrow k_{n+1}(E)_\S/i_{E/F}(k_{n+1}(F))
 $$
for which $\tr_{E/F}\circ\phi=\id$.
 This is done as follows.

 By our assumptions, the map $\tr_{E/F}\:k_n(E)_\S/i_{E/F}k_n(F)
\rarrow k_n(F)$ is an isomorphism.
 Therefore, we can define a homomorphism
 $$
  \Phi\:k_1(F)\ot k_n(F)\lrarrow k_{n+1}(E)_\S/i_{E/F}(k_{n+1}(F))
 $$
by the rule $\Phi(b\ot\alpha)=\overline{b\beta}$, where
$\beta\in k_n(E)$ is such that $\tr_{E/F}\beta=\alpha$.
 Let us verify that $\Phi$ can be factorized through $k_{n+1}(F)$.
 It suffices to show that $\Phi((1-a)\ot\alpha)=0$ whenever
$\alpha\in k_n(F)$~is divisible by the class $(a)\in k_1(F)$ of
an element $a\in F\setminus\{0,1\}$.

 Let us first assume that $a$ is not an $\l$\+th power in~$E$.
 Consider the field $K=F(\sqrt[\l]{a})$; let $L=EK$ be
the field generated by $E$ and~$K$ over~$F$.
 Then we have $\tr_{L/K}i_{L/E}\beta=i_{K/F}\tr_{E/F}\beta=0$,
since $\tr_{E/F}\beta=\alpha$ is divisible by~$(a)$.
 By assumption, it follows that
$i_{L/E}\beta\.\in\. i_{L/K}k_n(K) + (1-\sigma)k_n(L)$,
where $\sigma$~is a generator of the group $\Gal(L/K)=\Gal(E/F)$.
 Now we have
\begin{multline*}
 (1-a)\cdot\beta \.=\. \tr_{L/E}((1-\sqrt[\l]{a}))\cdot \beta
 \.=\. \tr_{L/E}((1-\sqrt[\l]{a}) \cdot i_{L/E}\beta) \\
 \sub\. \tr_{L/E}(i_{L/K}k_{n+1}(K) + (1-\sigma) k_{n+1}(L))
 \.\sub\. i_{E/F}k_{n+1}(F) + (1-\sigma) k_{n+1}(E),
\end{multline*}
hence $\Phi((1-a)\ot\alpha)=0$.

 It remains to consider the case $\sqrt[\l]{a}\in E$.
 In this case we have
$\sum_{s=1}^\l\sigma^s(\beta)=i_{E/F}\tr_{E/F}\beta=i_{E/F}\alpha=0$
and therefore
 $$
  \textstyle
  (1-a)\cdot\beta \.=\. \sum_{s=1}^\l\sigma^s(1-\sqrt[\l]{a})\cdot
  \beta \.=\. \sum_{s=1}^\l(1-\sqrt[\l]{a})\cdot\sigma^{-s}(\beta)=0
 $$
modulo $(1-\sigma)k_{n+1}(E)$.

 Thus we have constructed the desired map $\phi$.
 Since the field~$F$ has no finite extensions of degree
less than~$\l$, we have $k_{n+1}(E)=k_1(F)\cdot k_n(E)$
(see~\cite{BT} and section~4); hence the map~$\phi$ is surjective.
 It is obvious that $\tr_{E/F}\circ\phi=\id$.
\end{proof}

 In the remaining part of this section (which is included for
the sake of completeness of the exposition) we closely
follow~\cite[Section~5]{Voev}.

\begin{corollary}
 \.Assume that for any finite field extensions $F\sub K\sub L$
with $[L:K]=\l$ the following sequence of Milnor K\+theory
groups is exact
 $$
  k_n(K)\overset{i}\lrarrow k_n(L)_\S
  \overset{\tr}\lrarrow k_n(K)
  \overset{\!u\,\cdot\.}\lrarrow k_{n+1}(K),
 $$
where $u\in k_1(K)$ is the element corresponding to the cyclic
extension~$K\sub L$ and $\S=\Gal(L/K)$.
 Further assume that $k_{n+1}(F)=0$.
 Then one has $k_{n+1}(E)=0$ for
any finite field extension~$F\sub E$.
\end{corollary}

\begin{proof}
 Proceeding by induction, we only have to consider the case
when $[E:F]=\l$.
 Since $k_{n+1}(F)=0$, it follows from our exact sequence that
the transfer homomorphism $\tr_{E/F}\:k_n(E)\rarrow k_n(F)$
is surjective.
 Therefore, the conditions of Lemma~3 are satisfied and
we can conclude that $k_{n+1}(E)_\Sigma=0$.
 Thus $k_{n+1}(E)=(1-\sigma)k_{n+1}(E)=(1-\sigma)^\l k_{n+1}(E)=0$.
\end{proof}

\begin{proof}[Proof of Theorem 1]
 Let $F\sub K\sub L$ be finite field extensions with $[L:K]=\l$.
 We will apply Proposition~2 for $m=\l$, the pro-finite group
$G=G_K$, the subgroup $H=G_L$, the modules $T_{\l^2}=\mu_{\l^2}^
{\ot n}$ and $T_\l\simeq \Z/\l$ over~$G$, and the degree~$n$.
 It follows from the assumptions of Theorem~1 that the required
Bockstein homomorphisms vanish.
 From Proposition~2 we conclude that the conditions of Corollary~4
are satisfied for the field~$F$.
 Thus we have $k_{n+1}(E)=0$ for any field~$E$ finite over~$F$.

 Now assume that $\alpha$~is a nonzero element
of $H^{n+1}(G_F,\.\Z/\l)$.
 It is clear that there exist finite field extensions
$F\sub K\sub L$ such that $\res_{K/F}\alpha\ne0$, but
$\res_{L/F}(\alpha)=0$ and $[L:K]=\l$.
 Using Proposition~2 again, we see that
$\res_{K/F}\alpha=u_{L/K}\cup\beta$
for some $\beta\in H^n(G_K,\.\Z/\l)$.
 Since $k_n(K)\simeq H^n(G_K,\.\Z/\l)$ and $k_{n+1}(K)=0$,
it follows that $\res_{K/F}\alpha=0$.
 This contradiction proves that $H^{n+1}(G_F,\.\Z/\l)=0$.
\end{proof}

\Section{The Six-Term Cohomological Exact Sequence}

 In this section we develop a rather general setting for
exact sequences of cohomology under the assumptions of vanishing
of Bockstein homomorphisms.
 The results below are actually valid for an arbitrary
cohomological functor on an abelian category endowed with
a central element with zero square, etc.
 We restrict ourselves to the cohomology of a pro-finite group
for the sake of simplicity of exposition only.

 Let us fix a pro-finite group~$G$ and a number~$m\ge 2$.
 We will use the following notation: $A_2$, \.$B_2$,~etc.\ will
denote free modules over $\Z/m^2$ endowed with a discrete action
of~$G$ and $A_1$, \.$B_1$,~etc.\ will be the corresponding
$G$\+modules over~$\Z/m$, defined as $A_1=m A_2=A_2/m A_2$.
 For a $G$\+module~$M$, we will simply write $H^n(M)=H^n(G,\.M)$.
 The Bockstein homomorphisms corresponding to exact triples
 $$
  0\lrarrow M_1\overset{\tau}\lrarrow
  M_2\overset{\pi}\lrarrow M_1\lrarrow 0,
  \qquad \tau\pi=m
 $$
will be denoted by $\beta^{\.n}_M\:H^n(M_1)\rarrow H^{n+1}(M_1)$.

\begin{lemma}
 \.Let $0\to X_2\overset{i}\rarrow Y_2\overset{p}\rarrow Z_2
\to 0$ be an exact triple of $G$\+modules over~$\Z/m^2$,
$\.0\to X_1\rarrow Y_1\rarrow Z_1\to0$ be the corresponding
exact triple of $G$\+modules over~$\Z/m$, and
$\d_{XZ}\: H^n(Z_1)\rarrow H^{n+1}(X_1)$ be the induced
homomorphism of pro-finite group cohomology.

 Let $h$ be a homotopy map on the complex $X_2\rarrow Y_2\rarrow
Z_2$ such that $dh+hd=m\cdot\id$, where $d=(i,\.p)$.
 Then the induced map $h_1$ is an endomorphism of degree~$-1$
of the complex $X_1\rarrow Y_1\rarrow Z_1$.
 It follows that there is a map $h_{XZ}\:Z_1\rarrow X_1$
such that $h_{XZ}\circ p_1=h_1\:Y_1\rarrow X_1$ and
$i_1\circ h_{XZ}=-h_1\:Z_1\rarrow Y_1$.

 Then the following equality holds:
$\d_{XZ}=\beta_X\circ H(h_{XZ}) - H(h_{XZ})\circ\beta_Z$.
\end{lemma}

\begin{proof}
 The desired equality of homomorphisms of cohomology is to be
understood as the corollary to an identity in the $G$\+module
extension group $\Ext^1_{\Z/m^2[G]}(Z_1,\.X_1)$, where
$\d_{XZ}$ and the betas represent certain extension classes
and $h_{XZ}$~is just a map of modules.
 This identity in the group of extension classes is verified
as follows.

 The composition of the Bockstein extension $X_1\rarrow X_2
\rarrow X_1$ with $h_{XZ}\:Z_1\rarrow X_1$ is an extension
of the form $X_1\rarrow X_2\sqcap_{X_1}Z_1\rarrow Z_1$, while
the composition of $Z_1\rarrow Z_2\rarrow Z_1$ with~$h_{XZ}$
is $X_1\rarrow X_1\sqcup_{Z_1}Z_2\rarrow Z_1$
(where $\sqcap$ and~$\sqcup$ denote the relative product
and coproduct, respectively). 
 We have to check that the difference of these extensions is
isomorphic to $X_1\rarrow Y_1\rarrow Z_1$.
 The middle term of the desired difference can be defined as
the homology module of the sequence
$$ X_1 \overset{(\tau,0,\id,0)}\longrarrow
(X_2\sqcap_{X_1}Z_1) \op (X_1\sqcup_{Z_1}Z_2) 
\overset{(0,\id,0,-\pi)}\lolongrarrow Z_1.
$$
Then the arrows $j$ and $q$ forming this extension are induced
by
$$X_1 \overset{(\tau,0,0,0)}\lorarrow
(X_2\sqcap_{X_1}Z_1) \op (X_1\sqcup_{Z_1}Z_2)
\quad \text{and} \quad
(X_2\sqcap_{X_1}Z_1) \op (X_1\sqcup_{Z_1}Z_2) 
\overset{(0,0,0,\pi)}\longrarrow Z_1,$$
respectively.
 Let us construct an homomorphism~$f$ from $Y_1$ to the
homology module in the middle.
 Given an element $y_1\in Y_1$, we choose its preimage
$y_2\in Y_2$ and set $f(y_1)=(h(y_2),\.p_1(y_1))\op
(h_1(y_1),\.p(y_2))$.
 It is obvious that $(0,\id,0,-\pi)\circ f(y_1)=
p_1(y_1)-\pi p(y_2)=0$, since $p$ commutes with $\pi$.
 Now let us check that $f$ is well-defined.
 It suffices to restrict oneself to the case
when $y_1=0$ and $y_2=\tau(y'_1)$.
 By the above formula, we have $f(y_1)=(h\tau(y'_1),\.0)\op
(0,\.p\tau(y'_1))=(\tau h_1(y'_1),\.0)\op (0,\.\tau p_1(y'_1))=
(\tau h_1(y'_1),\.0)\op (h_{XZ}p_1(y'_1),\.0)=
(\tau h_1(y'_1),\.0)\op (h_1(y'_1),\.0)=
(\tau,0,\id,0)(h_1(y'_1))$.
 It remains to check commutativity of the triangles formed
by the maps $f$, $i_1$, $j$ and $f$, $p_1$, $q$.
 For any $x_1\in X_1$, choose a preimage $x_2\in X_2$ and
take $y_2=i(x_2)$ for $y_1=i_1(x_1)$.
 By the same formula, we have $f(i_1(x_1))=(hi(x_2),\.0)
\op (h_1i_1(x_1),\.0)=(mx_2,\.0)\op (0,0)=(\tau(x_1),\.0)
\op (0,0) = j(x_1)$.
 Finally, for any $y_1\in Y_1$ with a preimage $y_2\in Y_2$,
we have $(q\circ f)(y_1)=\pi p(y_2)=p_1(y_1)$.
\end{proof}

\begin{theorem}
 \.Let $G$ be a pro-finite group, $m\ge2$ and $n\ge0$ be some
integers, and
 $$
  0\lrarrow A_2\lrarrow B_2\lrarrow C_2\lrarrow D_2\lrarrow 0
 $$
be a 4\+term exact sequence of free $\Z/m^2$\+modules with
a discrete action of\/~$G$ in it.
 Suppose that we are given a homotopy map~$h$ in this
exact sequence such that $dh+hd=m$.
 Assume that the Bockstein maps
 $$
  \beta_X^n\: H^n(G,X_1)\lrarrow H^{n+1}(G,X_1),
  \qquad X_1=mX_2=X_2/m
 $$
vanish for all $X=A$, $B$, $C$, or~$D$ and the given~$n$.
 Then there is a 6\+term exact sequence of the form
 \begin{multline*}
   H^n(B_1\op D_1) \overset{(d_1,\,h_1)}\longrarrow
   H^n(C_1) \overset{d_1}\lrarrow H^n(D_1)        \\
   \overset{\nu}\lrarrow H^{n+1}(A_1)
   \overset{d_1}\lrarrow H^{n+1}(B_1)
   \overset{(h_1,\,d_1)}\longrarrow H^{n+1}(A_1\oplus C_1).
 \end{multline*}
The middle arrow~$\nu$ comes from a certain extension
 $$
  0\lrarrow A_1\lrarrow N\lrarrow D_1\lrarrow 0
 $$
of discrete $G$\+modules over~$\Z/m$, whose compositions with
the morphisms $A_1\rarrow B_1$ and $C_1\rarrow D_1$ are trivial.
 If the group~$G$ acts on the original exact quadruple of modules
through its finite quotient group~$\Gamma$, then $N$ is
a $\Z/m[\Gamma]$\+module.
\end{theorem}

\begin{proof}
 Let $L_2$ denote the image of the differential $B_2\rarrow C_2$.
 It is clear that $L_2$ is a free module over $\Z/m^2$.
 Now our 4\+term exact sequence of $G$\+modules consists of two
exact triples: $A_2\rarrow B_2\rarrow L_2$ and $L_2\rarrow C_2
\rarrow D_2$.
 Reducing modulo~$m$, we see that the map~$h_1\:C_1\rarrow B_1$
defines a homomorphism from the exact triple
$L_1\rarrow C_1\rarrow D_1$ to $A_1\rarrow B_1\rarrow L_1$.
 The corresponding ``intermediate'' induced extension provides
the desired exact triple $A_1\rarrow N\rarrow D_1$.

 Now the 6\+term sequence is constructed; it remains to check
that it's exact.
 Let us introduce notation for maps $p\:B_2\rarrow L_2$ and
$i\:L_2\rarrow C_2$; the same maps reduced modulo~$m$ will be
be denoted by $p_1$ and~$i_1$.
 It is easy to see that the homotopy map $h$ induces analogous
homotopies on the exact triples $A_2\rarrow B_2\rarrow L_2$
and $L_2\rarrow C_2 \rarrow D_2$.
 Therefore, both triples satisfy the conditions of Lemma~5.
 We have connecting maps
$\d_{AL}\:H^n(L_1)\rarrow H^{n+1}(A_1)$ and
$\d_{LD}\:H^n(D_1)\rarrow H^{n+1}(L_1)$
in the pro-finite group cohomology;
the cohomology maps induced by the module maps
$h_{AL}\:L_1\rarrow A_1$ and $h_{LD}\:D_1\rarrow L_1$
will be denoted simply by $h_{AL}$ and $h_{LD}$.
 According to Lemma~5, there are two identities:
$\d_{LD}=\beta_L \circ h_{LD} - h_{LD} \circ \beta_D$ and
$\d_{AL}=\beta_A \circ h_{AL} - h_{AL} \circ \beta_L$.
 We also know that $\nu = \d_{AL}\circ h_{LD} =
- h_{AL}\circ \d_{LD}$.

 Assume that $x\in H^n(D_1)$ and $\nu(x)=0$.
 Let us prove that $x\in d_1H^n(C_1)$.
 We have $\d_{AL}h_{LD}(x)=0$, hence $h_{LD}(x)=p_1(z)$
for some $z\in H^n(B_1)$.
 Since $p_1$ is induced by a map of modules over $\Z/m^2$,
it commutes with betas.
 So it follows from our assumptions about vanishing of
Bockstein homomorphisms that
$0=p_1\beta_B(z)=\beta_L p_1(z)=\beta_L h_{LD}(x) =
\beta_L h_{LD}(x) - h_{LD}\beta_D(x) = \d_{LD}(x)$.
 Thus $x\in d_1H^n(C_1)$.

 Next assume that $w\in H^{n+1}(A_1)$ and $d_1(w)=0$.
 Then $w=\d_{AL}(u)$ for some $u\in H^n(L_1)$, hence
$w=\beta_A h_{AL}(u)-h_{AL}\beta_L(u)=-h_{AL}\beta_L(u)$.
 The map~$i_1$ commutes with betas for the same reason
as above, so according to our assumption about Bockstein
maps we have $i_1\beta_L(u)=\beta_C i_1(u)=0$.
 It follows that $\beta_L(u)=\d_{LD}(x)$ for some
$x\in H^n(D_1)$.
 Thus $w=-h_{AL}\beta_L(u)=-h_{AL}\d_{LD}(x)=\nu(x)$.

 Now assume that $y\in H^n(C_1)$ and $d_1(y)=0$.
 Let us prove that $y\in d_1H^n(B_1) + h_1H^n(D_1)$.
 Since $d_1(y)=0$, we have $y=i_1(u)$ with $u\in H^n(L_1)$.
 Then $i_1\beta_L(u) = \beta_C i_1(u)=0$, hence
$\beta_L(u)=\d_{LD}(x)$ for some $x\in H^n(D_1)$.
 Therefore, $\d_{AL}(u)=\beta_A h_{AL}(u) -
h_{AL}\beta_L(u) = -h_{AL}\beta_L(u) = -h_{AL}\d_{LD}(x)
= \d_{AL}h_{LD}(x)$.
 It follows that $\d_{AL}(u-h_{LD}(x))=0$ and
$u - h_{LD}(x) = p_1(z)$ for some $z\in H^n(B_1)$.
 Applying $i_1$, we obtain $y = i_1p_1(z)+i_1h_{LD}(x) =
d_1(z) - h_1(x)$.

 Finally, assume that $t\in H^{n+1}(B_1)$ is such that
$d_1(t)=0$ and $h_1(t)=0$.
 Then $i_1p_1(t)=d_1(t)=0$, hence $p_1(t)=\d_{LD}(x)$
for some $x\in H^n(D_1)$.
 We have $\d_{AL}h_{LD}(x)=-h_{AL}\d_{LD}(x)=
-h_{AL}p_1(t)=-h_1(t)=0$, so $h_{LD}(x)=p_1(z)$ with
$z\in H^n(B_1)$.
 It follows that $p_1(t)=\d_{LD}(x)=\beta_Lh_{LD}(x)-
h_{LD}\beta_D(x)=\beta_Lh_{LD}(x)=\beta_Lp_1(z)=
p_1\beta_B(z)=0$ and therefore $t\in d_1H^{n+1}(A_1)$.
\end{proof}

\begin{remark}
 The middle arrow~$\nu$ does not depend on the choice
of a homotopy map~$h$, provided that $h$ satisfies
the conditions of Theorem~6.
 Indeed, let us consider another homotopy map $h'=h+t$,
where $dt+td=0$.
 Then the induced maps $t_{AL}$ and $t_{LD}$ commute
with Bockstein homomorphisms, and according to the above
computations we have $\nu'-\nu = \d_{AL}\circ t_{LD} =
\beta_A \circ h_{AL}\circ t_{LD} - h_{AL} \circ \beta_L
\circ t_{LD} = \beta_A \circ h_{AL} \circ t_{LD} -
h_{AL} \circ t_{LD} \circ \beta_D = 0$ on $H^n(D_1)$.
 If $A_2$ and $D_2$ are permutational modules (as in
all applications below), then $h_{AL} \circ t_{LD}$
can be lifted to a map $D_2\rarrow A_2$ and it follows
that the extension $A_1\rarrow N\rarrow D_1$ does not
depend on~$h$ either.
\end{remark}

\Section{The Relative Conjecture of Bloch and Kato}

 In this section we will prove two theorems providing exact
sequences for Galois cohomology of cyclic field extensions,
both of them generalizations of Proposition~2.
 The first one will be essentially the ``relative conjecture''
of Bloch and Kato~\cite[Conjecture~3.1]{BK}.
 In particular, the proof of Theorem~1 will be completed.

 Throughout this section, we will use the following notation:
$m\ge2$ is an integer, $G$~is a pro-finite group,
$T_{m^2}$ is a free module over $\Z/m^2$ endowed with
an action of~$G$, and $T_m=m T_{m^2}=T_{m^2}/m T_{m^2}$
is the corresponding $G$\+module over~$\Z/m$.
 Furthermore, $n\ge 0$ is another integer, and it is assumed
that the Bockstein homomorphisms
 $$
  H^n(H,\.T_m)\lrarrow H^{n+1}(H,\.T_m), \qquad H\sub G
 $$
vanish for all relevant open subgroups~$H$ in~$G$, including
$G$~itself.
 In applications to Galois cohomology it will be always meant
that $G=G_F$ is the absolute Galois group of a field~$F$ and
$T_{m^2}=\mu_{m^2}^{\otimes n}$ is the cyclotomic module.

\begin{proposition}
 \.Let $H\sub G$ be a normal subgroup such that $\S=G/H$
is a cyclic group of an order~$k$ dividing~$m$.
Then the following sequence of cohomology groups is exact
 \begin{multline*}
   H^n(G,\.T_m) \overset{m/k\cdot\res}\longrarrow
   H^n(H,\.T_m)_\S \overset{\cor}\lrarrow H^n(G,\.T_m) \\
   \overset{\!u\,\cup\.}\lrarrow H^{n+1}(G,\.T_m)
   \overset{\res}\lrarrow H^{n+1}(H,\.T_m)^\S
   \overset{m/k\cdot\cor}\longrarrow H^{n+1}(G,\.T_m),
 \end{multline*}
where $u\in H^1(G,\.\Z/m)$ is the class corresponding
to the character $G\rarrow\S\hookrightarrow \Z/m$.
\end{proposition}

\begin{proof}
 Consider the following very well known 4\+term exact sequence
of modules over the cyclic group~$\S$
 $$
  0\lrarrow\Z \overset{1+\dsb+\sigma^{k-1}}\longlongrarrow
  \Z[\S] \overset{1-\sigma}\lorarrow \Z[\S]
  \overset{\sigma^i\.\longmapsto 1}\longrarrow \Z \lrarrow0,
 $$
where $\sigma$ is a generator of~$\Sigma$.
 Let us construct a homotopy map~$h_k$ such that $dh_k+h_kd=k$.
 Set $f(x)$ to be the polynomial for which $k-(1+\dsb+x^{k-1})
=(1-x)f(x)$.
 The leftmost component of~$h_k$ sends $\sigma^i$ to~$1$
for all~$i$, the middle part is given by the multiplication
with~$f(\sigma)$ in the group ring~$\Z[\S]$, and the rightmost
component sends~$1$ to $1+\ds+\sigma^{k-1}$.
 For the purposes of our proof it suffices to put $h=m/k\cdot h_k$,
tensor the above exact quadruple with~$T_{m^2}$, and apply Theorem~6.
\end{proof}

 The proof of Theorem~1 is now finished.

\begin{proposition}
 \.Let $H\sub G$ be a normal subgroup such that $\S=G/H$
is a cyclic group of an order~$k$ divisible by~$m$.
 Let $\Pi\sub\Sigma$ be the cyclic subgroup of order~$m$ and
$K\sub G$ be kernel of the projection $G\rarrow \S/\Pi$.
 Then the following sequences of cohomology groups are exact
 $$
   H^n(H)_\S \overset{\cor}\lrarrow H^n(K)_{\S/\Pi}
   \overset{\cor{}\circ\,\.{u\.\cup}}\longrarrow H^{n+1}(G)
   \overset{\res}\lrarrow {}_{\!}H^{n+1}(H)^\S
   \overset{\cor}\lrarrow H^{n+1}(K)^{\S/\Pi}    
 $$
and
 $$
   H^n(K)_{\S/\Pi} \overset{\res}\lrarrow H^n(H)_\S
   \overset{\cor}\lrarrow H^n(G)
   \overset{\!{u\.\cup}\,\,\circ{}\res}\longrarrow
   H^{n+1}(K)^{\S/\Pi}\overset{\res}\lrarrow H^{n+1}(H)^\S,
 $$
where we write $H^i(I)=H^i(I,\.T_m)$ for any subgroup $I\sub G$,
and $u\in H^1(K,\.\Z/m)$ is the cohomology class corresponding
to the character $K\rarrow\Pi\simeq\Z/m$.
\end{proposition}

\begin{proof}
 The first of the desired exact sequences holds due to
the following exact quadruple~$Q_\S$ of permutational
$\S$\+modules
 $$
  \Z \lrarrow \Z[\S]\op \Z \lrarrow \Z[\S]\op \Z[\S/\Pi]
  \lrarrow \Z[\S/\Pi],
 $$
where the first map is $(1+\dsb+\sigma^{k-1}, -m)$,
the third map is $(\pr_\Pi, -1+\bar\sigma)$, and the middle
arrow is given by the matrix~$U$ with components
$U_{11}=1-\sigma$, $U_{12}=0$, $U_{21}=\pr_\Pi$, and
$U_{22}=1+\dsb+{\bar\sigma}^{k/m-1}$.
 Here $\sigma$~is a generator of~$\S$, by $\bar\sigma$
we denote its image in~$\S/\Pi$, and $\pr_\Pi$~is
the projection map $\Z[\S]\rarrow \Z[\S/\Pi]$.
 Let us construct a homotopy map~$h$ with $dh+hd=m$.
 Set $f(x)$ to be the polynomial for which 
$m-(1+x^{k/m}+\dsb+x^{(m-1)k/m})=(1-x)f(x)$.
 Denote by $\#\Pi$ the lifting map
$\Z[\S/\Pi]\rarrow \Z[\S]$ which sends $\bar\sigma^i$ to
$(1+\sigma^{k/m}+\dsb+\sigma^{(m-1)k/m)})\sigma^i$.
 Then the leftmost piece of $h$ is $(0,-1)$, the rightmost
piece is $(\#\Pi,-f(\bar\sigma))$, and the middle map
is given by the matrix $H$ with entries $H_{11}=f(\sigma)$,
$H_{12}=\#\Pi$, $H_{21}=0$, and $H_{22}=0$.
  Using this exact quadruple, one argues as in the proof
of Proposition~7 above.
 To obtain the second exact sequence of cohomology, one
has to use the dual exact quadruple of $\S$\+modules
$\Hom_\Z(Q_\S,\.\Z)$.
\end{proof}

\Section{The Bass--Tate Lemma and Applications to Annihilator Ideals}

 The next theorem provides a more sophisticated version of
the classical technique known as the ``Bass--Tate lemma''~\cite{BT}.

\begin{theorem}
 \.Let $E/F$ be a separable extension of fields of degree~$k$.
 Assume that the field~$F$ has no nontrivial finite extensions
of degrees less or equal to~$k/2$.
 Then $\KM_{\ge1}(E)=\bigoplus_{i=1}^\infty\KM_n(E)$
is a \textit{quadratic module} over the ring $\KM_*(F)$,
i.~e., a graded module generated in degree~$1$ with
relations in degree~$2$.
\end{theorem}

\begin{lemma}
 \.Let $E=F[t]$ be a finite extension of fields of degree~$k$,
where $t\in E$ is a generator.
 Then any element of the field~$E$ can be represented as
a fraction of polynomials $f(t)/g(t)$, where $f$ and~$g$
are two polynomials of degree~$\le k/2$ each.
\end{lemma}

\begin{proof}
 Let $V\sub E$ be the vector subspace of all elements of
the form~$f(t)$, where $\deg f\le k/2$.
 Then $2\dim V\ge k+1$, thus for any element $e\in E$ one has
$V\cap eV\ne0$.
 This simple argument was communicated to the author by
Alexander Vishik.
\end{proof}

\begin{lemma}
 \.Let $A\rarrow B$ be a surjective homomorphism of abelian
groups.
 Then the kernel $I^*$ of the induced map of exterior rings
$\Lambda^*(A)\rarrow \Lambda^*(B)$ is generated in degree~1
as an ideal in $\Lambda^*(A)$.
\end{lemma}

\begin{proof}
 Consider an exterior expression of elements of~$A$ which
becomes zero in $\Lambda^*(B)$.
 There is a chain of transformations (by the rules of exterior
algebra) of exterior expressions of elements of~$B$ connecting
this image with the zero.
 It is clear that one would be able to lift this chain of
transformations to the original expression in~$A$ if one
only were allowed to freely change the elements of $A$ in
the expressions with other elements of $A$ with the same
image in~$B$.
 In other words, the ideal $I^*$ is additively generated by
exterior monomials $a_1\wedge\dsb\wedge a_n$ where one of
the $a_i$'s belongs to the kernel of $A\rarrow B$.
 But this is exactly what we wanted to prove.
\end{proof}

\begin{proof}[Proof of Theorem 9]
 Consider the field of rational functions~$F(x)$ in
a transcendental variable~$x$.
 Let $D\sub F(x)^*$ be the multiplicative subgroup generated
by polynomials of degree~$\le 1$.
 Let $R_*=\bigoplus_{n=0}^\infty R_n$ be the graded
anti-commutative ring generated by $R_1=D$ with relations
$\{f,\.1-f\}=0$ for all $f$, $1-f\in D$ and $\{f,-f\}=0$
for all $f\in D$.
 Arguing as in~\cite[Theorem~I.5.1]{BT}, it is not hard to
show that $R_*$~is a free $\KM_*(F)$\+module with a basis
consisting of~$1$ and $\{x-a\}$ for all $a\in F$.
 Hence the $\KM_*(F)$\+module $R_{\ge1}$ is the direct sum
of $\KM_{\ge1}(F)$ and a free module, so it is quadratic.

 Now let us choose an element $t\in E$ such that $E=F[t]$.
 Let $\pi(x)\in F[x]$ be the irreducible equation of
the element~$t$.
 Obviously, there is a homomorphism $p_\pi\:R_*\rarrow \KM_*(E)$
assigning $\{f \bmod \pi\}$ to~$\{f\}$ for any $f\in D$.
 According to Lemma~10, this homomorphism is surjective.
 Moreover, the ideal $J_*=\ker p_\pi\sub R_*$ is generated
by its first-degree component~$J_1$.
 Indeed, the kernel of the map of exterior algebras
$\Lambda^*(D)\rarrow \Lambda^*(E^*)$ is generated in degree~1
by Lemma~11, and it remains to show that any Steinberg element
in $\Lambda^2(E^*)$ can be lifted to one of the relations
defining $R$.
 This follows from the fact, again established by Lemma~10,
that for any $e\in E\setminus\{0,1\}$ there exists
$f\in D\sub F(x)^*$ such that $1-f\in D$
and $\pi(f)=e$.

 Since $R_*$~is anti-commutative, we have
 $$
  J_*=R_*\cdot J_1 = (\KM_*(F) + \KM_*(F)\cdot R_1)\cdot J_1
  = \KM_*(F) \cdot (J_1 + R_1 J_1).
 $$
 Therefore, the $\KM_*(F)$\+module $J_*$ is generated by $J_1$
and~$J_2$, while the module $R_{\ge1}$ is quadratic.
 Thus the quotient module $\KM_{\ge1}(E)=R_{\ge1}/J_*$
is also quadratic.
\end{proof}

 For reader's convenience, the statement of Bass--Tate lemma
proper is formulated here, together with an improvement by
Becher~\cite{Bech}.

\begin{proposition}
 \.Let $E/F$ be a separable field extension of degree~$k=\l^s$,
where $\l$~is a prime number.
 Assume that either (a)~$\l=2$, or (b)~the field~$F$ has no finite
extensions of degree different from powers of\/~$\l$ and
less or equal to~$k/2$, or (c)~$E/F$ is a tower of extensions
of degree~$\l$ such that the largest proper subfield of\/~$E$ in
this tower has no nontrivial finite extensions of degree~$\le \l/2$.
 Then the Milnor ring~$\KM_*(E)$ is generated in degree~$\le s$
as a module over~$\KM_*(F)$.
\end{proposition}

\begin{proof}
 The case (a) is covered by the result of~\cite{Bech}.
 The case (b) follows from Lemma~10 and~\cite[Lemma~I.5.2]{BT}.
 In the case~(c) one argues by induction in~$s$, applying~(b)
to extensions of degree~$\l$.
\end{proof}

 From now on and for the rest of this section we will assume
that the Milnor--Kato conjecture holds for all relevant fields,
primes, and degrees.
 Under this assumption we will deduce corollaries from
Proposition~7, Proposition~8, and Theorem~9.

\begin{corollary}
 \.Let $\l$ be a prime number, $m$ be a power of\/~$\l$, and $F$~be
a field containing a primitive root of unity of degree~$m$.
 Consider the ring $H^*(G_F,\.\Z/m)$.
 Let $u\in H^1(G_F,\.\Z/m)$ be an element of additive order~$k=\l^s$
and $E/F$ be the corresponding cyclic field extension of degree~$k$.
 Assume that the maximal proper subfield of\/~$E$ over~$F$ has
no nontrivial finite extensions of degree $\le \l/2$.
 Then the ideal $\Ann(u)$ in the ring $H^*(G_F,\.\Z/m)$
is generated by elements of degree less or equal to~$s$.
 In particular, if\/ $m=\l$ is a prime and $F$ has no nontrivial
extensions of degree~$\le \l/2$, then $\Ann(u)$ is generated
in degree~$1$.
\end{corollary}

\begin{proof}
 According to Proposition~7, we have an exact sequence of
$H^*(G_F,\.\Z/m)$-modules
 $$
  H^*(G_E,\.\Z/m) \overset{\cor}\lrarrow H^*(G_F,\.\Z/m)
  \overset{\!u\,\cup\.}\lrarrow H^*(G_F,\.\Z/m).
 $$
 But $H^*(G_E,\.\Z/m)$ is generated in degree~$\le s$ as
a module over $H^*(G_F,\.\Z/m)$; this follows from
the Milnor--Kato conjecture and Proposition~12.
\end{proof}

\begin{corollary}
 \.Let $\l$ be a prime number, $m$ and $k$~ be powers of\/~$\l$,
and $F$~be a field containing a primitive root of unity of
degree~$m$.
 Let $L/F$ be a cyclic extention of degree~$k$.
 Consider the kernel~$J$ of the ring homomorphism
$H^*(G_F,\.\Z/m) \rarrow H^*(G_L,\.\Z/m)$.
 If $k$ is divisible by $m$, consider the intermediate extension
$F\sub E\sub L$ of degree~$k/m$ over $F$ and  assume that
the maximal proper subfield of\/~$E$ over~$F$ has
no nontrivial finite extensions of degree $\le \l/2$.
 Then the ideal~$J$ in the ring $H^*(G_F,\.\Z/m)$
is generated by elements of degree less or equal to $s+1$,
where $k/m=\l^s$ if\/ $k$~is divisible by~$m$ and $s=0$ if
$m$~is divisible by~$k$.
\end{corollary}

\begin{proof}
 The case when $k$~is divisible by~$m$ follows from Proposition~8
in the way analogous to the proof of Corollary~13.
 Namely, one has to use the exact sequence
 $$
   H^n(G_E,\.\Z/m) \overset{\cor{}\circ\,\.{u\.\cup}}
   \lolongrarrow H^{n+1}(G_F,\.\Z/m)
   \overset{\res}\lrarrow H^{n+1}(G_L,\.\Z/m),
 $$
where $u\in H^1(G_E,\.\Z/m)$ is the class of the extension~$L/E$.
 The case when $m$~is divisible by~$k$ follows directly from
Proposition~7.
\end{proof}

 The next corollary is the main result of this section.

\begin{corollary}
 \.Let $\l$~be a prime number and $F$~be a field containing
a primitive root of unity of degree~$\l$.
 Assume that $F$~has no nontrivial finite extensions of
degree~$\le \l/2$ (note that for $\l=2$ or $\l=3$ this
always holds).
 Then for any element $u\in H^1(G_F,\.\Z/\l)$ the ideal
$\Ann(u)$ in the cohomology ring $H^*(G_F,\.\Z/\l)$
is a quadratic module over this ring, i.~e., it is
isomorphic to a graded module with generators and relations,
where generators are in degree~$1$ and relations in degree~$2$.
\end{corollary}

\begin{proof}
 Let $E/F$~be the cyclic extension of degree~$\l$ corresponding
to the class~$u$ and $\S$~be its Galois group.
 Consider the spectral sequence related to the maximal
($\l$\+term) filtration of the $G_F$\+module $\Z/\l\.[\S]$.
 It follows from the exact sequence
 $$
   H^n(G_E,\.\Z/\l) \overset{\cor} \lrarrow H^n(G_F,\.\Z/\l)
   \overset{\!u\,\cup\.} \lrarrow H^{n+1}(G_F,\.\Z/\l)
   \overset{\res}\lrarrow H^{n+1}(G_E,\.\Z/\l)
 $$
that this spectral sequence degenerates at the term~$E_2^{*,*}$.
 In other words, there is a filtration $V^1\supset\dsb\supset
V^\l$ on the $H^*(G_F,\.\Z/\l)$\+module $H^*(G_E,\.\Z/\l)$
with associated quotient modules of the form
$V^1/V^2 \simeq \ker({u\.\cup})$, $\.V^i/V^{i+1} \simeq
\ker({u\.\cup})/\im({u\.\cup})$ for $i=2$, \ds, $\l-1$,
and  $V^\l \simeq \coker({u\.\cup})$.
 From this filtration it follows that
the $H^*(G_F,\.\Z/\l)$\+modules $H^{\ge1}(G_E,\.\Z/\l)$
and $\Ann(u)$ are quadratic simultaneously.
 (They are also simultaneously Koszul, see~\cite{Pos}.)
 Now it remains to apply Theorem~9.
 Alternatively, one can argue as in the proof of Corollary~20
below.
\end{proof}

\begin{conjecture}
 \.Let $\l$~be a prime number and $F$~be a field containing
a primitive root of unity of degree~$\l$.
 Then for any element $u\in H^1(G_F,\.\Z/\l)$ the
$H^*(G_F,\.\Z/\l)$\+module $\Ann(u)\sub H^*(G_F,\.\Z/\l)$
is Koszul.
 For any cyclic extension~$E/F$ of degree~$\l$,
the $H^*(G_F,\.\Z/\l)$\+module $H^{\ge1}(G_E,\.\Z/\l)$
is Koszul (see~\cite{Pos} for the definition).
\end{conjecture}

\Section{Dihedral Field Extensions}

 Using Theorem~6 and assuming the Milnor--Kato conjecture, one can
construct some exact sequences of Galois cohomology for any finite
field extension whose Galois\- group admits a ``nontrivial enough''
exact quadruple of permutational representations over it.
 Unfortunately, examples of finite groups with such exact
quadruples of representations that I am aware of are very few.
 All of them turn out to be finite subgroups
of $\operatorname{PGL}_2({\mathbb C})$.
 These include cyclic groups, dihedral groups, the biquadratic group
$\Z/2\times\Z/2$, and the symmetric group ${\mathbb S}_4$.
 The goal of the next two sections is to construct those exact quadruples
for biquadratic and dihedral groups and deduce some corollaries, including
the conjectures of Merkurjev--Tignol~\cite{MT} and Kahn~\cite{Kahn}
about biquadratic field extensions.
 Throughout these sections we will assume that the characteristic of our
fields is not equal to~$2$ and the Milnor--Kato conjecture holds for $\l=2$.
 This is called the Milnor conjecture, and since it is proven by Voevodsky
already~\cite{Voev}, our results are in fact unconditional.

 We will sometimes use the notation $H^i(F)=H^i(G_F,\.\Z/2)$.

 Let $\Delta=\Delta_k$ be a dihedron of order~$k$, i.~e.,
a regular polygon with $k$~vertices and $k$~edges.
 We will denote by~$\Delta'$ the set of all vertices of~$\Delta$
and by~$\Delta''$ the set of all its edges.
 Let $\Gamma=\Gamma_k$~be the group of all automorphisms
(symmetries) of~$\Delta$ and $\S\sub\Gamma$ be the subgroup of
orientation-preserving automorphisms (rotations).
 Then $\S$~is a cyclic group of order~$k$ and a normal subgroup
of order~$2$ in~$\Gamma$.

 Assume that $k$~is divisible by~$4$.
 Let $\Pi$~be the (only) subgroup of order~$2$ in~$\S$.
 Then the following quadruple of permutational representations
of~$\Gamma$ is exact
 \begin{equation}
  \tag{$Q_\Gamma$}
  \Z \lrarrow \Z[\Delta']\op \Z \lrarrow \Z[\Delta'']
  \op \Z[\Delta'/\Pi] \lrarrow \Z[\Delta''/\Pi],
 \end{equation}
where the arrows are defined by the formulas below.
 Choose an arbitrary generator $\sigma\in\Sigma$; let $\bar\sigma$
be its image in~$\S/\Pi$.
 Further, let $\sigma^{\pm1/2}$ be maps 
$\Delta'\rarrow\Delta''\rarrow\Delta'$, commuting
with $\S$ and inverted by the action of $\Gamma\setminus\Sigma$,
whose squares are $\sigma^{\pm1}$, and let $\bar\sigma^{\pm1}$
be the induced maps
$\Delta'/\Pi\rarrow\Delta''/\Pi\rarrow\Delta'/\Pi$.
 Denote by~$\pr_\Pi$ the projection maps modulo~$\Pi$ and
by~$\#\.X$ the sum of all elements of a set~$X$.
 Then the first map in this sequence is $(\#\.\Delta', -2)$,
the third one is $(\pr_\Pi, {}-\bar\sigma^{-1/2}-\bar\sigma^{1/2})$,
and the middle arrow is given by the matrix~$U$ with components
$U_{11}=\sigma^{-1/2}+\sigma^{1/2}$, $U_{12}=\#\.\Delta''$,
$U_{21}=\pr_\Pi$, and $U_{22}=\#\.\Delta'/\Pi$.

 Moreover, there exists a chain homotopy~$h$ such that
$dh + hd = 2$.
 Of the three maps constituting~$h$, the leftmost one is given
by the pair~$(-\pr_1, {}-k/2-1)$, where $\pr_1\:\Z[X]\to\Z$
sends every $x\in X$ to~$1$.
 To define the remaining two maps, one needs some preparatory
work.
 Let $f(x)$ be a Laurent polynomial solving the equation
$(x^{-1/2}+x^{1/2})^2 \.f(x) = 1 + \sum_{i=-k/2+1}^{k/2-1} x^i$.
 Set $\psi=(\sigma^{-1/2}+\sigma^{1/2})f(\sigma)\in\sigma^{1/2}
\.\Z[\S]$ and $\bar\psi=\psi\bmod \Pi \in\bar\sigma^{1/2}\.
\Z[\S/\Pi]$.
 The rightmost piece of~$h$ is equal to~$(1+\pi,\,\.
\bar\sigma^{1/2}\.\#\.\S/\Pi \.-\. \bar\psi)$, where
$\pi$~is the only nontrivial element of~$\Pi$.
 The middle arrow is given by the matrix~$H$ with components
$H_{11}=\psi$, $H_{12}=1+\pi$, $H_{21}=-\pr_1$, and $H_{22}=0$.

 Further useful exact quadruples can be obtained by
restricting~$Q_\Gamma$ to a dihedral subgroup of index~$2$.
 The group~$\Gamma_{k/2}$ can be embedded into~$\Gamma_k$
in two ways; namely, under the first embedding the action of
$\Gamma_{k/2}$ on $\Delta'_k$ is transitive, while under
the second embedding the action of $\Gamma_{k/2}$ on $\Delta''_k$ is.
 These embeddings correspond to the two essentially different
ways of defining a regular polygon with $k/2$~vertices in
terms of the given polygon~$\Delta_k$: one can make edges of
the new polygon out of every second edge of the original one,
or one can take every second vertice of the original polygon
to be the new polygon's vertices.
 Taking the restriction of~$Q_\Gamma$ with respect to the first
of these emdeddings, we obtain an exact quadruple of
$\Gamma_{k/2}$-modules
 \begin{equation}
  \tag{$Q^+_\Gamma$}
  \Z \!\.\lrarrow\!\. \Z[\Delta^{\!\.+}]\op \Z \!\.\lrarrow\!\.
  \Z[\Delta'] \op \Z[\Delta''] \op \Z[\Delta^{\!\.+}\!\./\Pi]
  \!\.\lrarrow \!\.\Z[\Delta'/\Pi]\op\Z[\Delta''/\Pi],
 \end{equation}
where $\Delta'$ and $\Delta''$ denote the sets of all vertices
and edges of the regular polygon~$\Delta_{k/2}$, while
$\Delta^{\!\.+}$ denotes a principal homogeneous space for 
the dihedral group~$\Gamma_{k/2}$ which is defined as the set
of all pairs (an edge of~$\Delta_{k/2}$; one of its two vertices).

 Now let $M/F$ be an extension of fields with a dihedral Galois
group~$\Gamma$.
 Then the category of intermediate field extensions~$K/F$
embeddable into~$M$ is equivalent to the category of sets endowed
with a transitive action of~$\Gamma$.
 Let $L'$ and~$L''$ be the intermediate fields corresponding
to the $\Gamma$\+sets $\Delta'$ and~$\Delta''$, let
$K'$ and~$K''$ be the fields corresponding to $\Delta'/\Pi$
and $\Delta''/\Pi$, and $E\supset F$ be the normal extension
of~$F$ corresponding to the subgroup $\Pi\sub\Gamma$.

\begin{proposition}
 \.For any dihedral field extension~$M/F$ of degree $[M\:F]$
divisible by~$8$, there is an exact sequence of Galois cohomology
of the form
 \begin{multline*}
   H^n(L'')\op H^n(K')
   \overset{(\cor_{L''/K''},\;\.\cor_{E/K''}\circ\res_{E/K'})}
   \Longlongrarrow H^n(K'')
   \overset{\cor_{K''/F}\circ\;\.{u_{L''/K''}\.\cup}}
   \Longrarrow H^{n+1}(F) \\
   \overset{\res_{L'/F}}\longrarrow H^{n+1}(L')
   \overset{(\cor_{M/L''}\circ\res_{M/L'},\,\,\cor_{L'/K'})}
   \Longlongrarrow H^{n+1}(L'')\op H^{n+1}(K').
 \end{multline*}
For any dihedral extension~$M/F$ of degree divisible by~$4$,
there is an exact sequence 
 \begin{multline*}
   H^n(L')\op H^n(L'')\op H^n(E)\overset{\cor^4}\lorarrow
   H^n(K')\op H^n(K'')
   \overset{\cor_{K^{(i)}/F}\circ\;\.{u_{L^{(i)}/K^{(i)}}\.\cup}}
   \Lolongrarrow H^{n+1}(F) \\
   \overset{\res}\longrarrow H^{n+1}(M) \overset{\cor^3}
   \longrarrow H^{n+1}(L')\op H^{n+1}(L'')\op H^{n+1}(E).
 \end{multline*}
\end{proposition}

\begin{proof}
 The proposition follows from Theorem~6 applied to the above
exact quadruples $Q_\Gamma$ and~$Q^+_\Gamma$, respectively.
 There are also dual exact sequences of cohomology,
corresponding to the dual exact quadruples
$\Hom_\Z(Q_\Gamma,\.\Z)$ and $\Hom_\Z(Q^+_\Gamma,\.\Z)$,
but we will not write them down here.
\end{proof}

\begin{remark}
 There is yet another exact quadruple of permutational
representations of a dihedral group~$\Gamma$, whose restrictions
to the subgroup of rotations and dihedral subgroups of index~$2$
are also of some interest.
 Let $\Delta_k$ be the dihedron of an arbitrary order~$k$ and
$\Pi\sub\S\sub\Gamma_k$ be a cyclic subgroup of order~$m$.
 Then there is a self-dual exact quadruple of the form
 $$
  \Z[\Delta'/\Pi] \lrarrow \Z[\Delta']\op \Z \lrarrow
  \Z[\Delta''] \op \Z \lrarrow \Z[\Delta''/\Pi]
 $$
if $k/m$ is even (and $\Delta''$ must be replaced
with~$\Delta'$ if $k/m$ is odd).
 The first map in this sequence is $(\#\.\Pi, -\pr_1)$,
the third map is $(\pr_\Pi, \.-\.\#\.\Delta''/\Pi)$,
and the middle arrow is given by the matrix~$U$ with components
$U_{11}=\sigma^{-(k/m-1)/2}+\dsb+\sigma^{(k/m-1)/2}$,
$U_{12}=\#\.\Delta''$, $U_{21}=\pr_1$, and $U_{22}=m$
in the above notation.
 It is not difficult to find a homotopy map~$h$
with $dh + hd = m$.
\end{remark}

\Section{Biquadratic Field Extensions}

 Let $\Theta=\{1,a,b,c\}$ denote an abelian group isomorphic
to~$\Z/2\times\Z/2$.
 We keep the notation of section~5 related to dihedral groups.
 The group~$\Theta$ can be embedded into the dihedral group
$\Gamma_4$ in two ways:
as above, under the first embedding the action of~$\Theta$
on~$\Delta'$ is transitive, while under the second embedding
the action of~$\Theta$ on~$\Delta''$ is.
 Taking the restrictions on~$\Theta$ of the exact
quadruple~$Q_\Gamma$ for~$\Gamma=\Gamma_4$, we obtain (after
one simple cancellation in the second case) two exact quadruples
of permutational representations of~$\Theta$.
 They have the form
\begin{equation}
 \tag{$Q_\Theta$}
  \Z \lrarrow \Z[\Theta]\op \Z \lrarrow \Z[\Theta/a] \op
   \Z[\Theta/b] \op \Z[\Theta/c] \lrarrow \Z^{\op3}/\Z
\end{equation}
and
 $$
  \Z \lrarrow \Z[\Theta/a]\op \Z[\Theta/b] \lrarrow \Z[\Theta]
  \op \Z \lrarrow \Z[\Theta/c],
 $$
where we denote for simplicity by~$a$ the subgroup~$\{1,a\}$ etc.
 All maps in these sequences can be defined directly in
a very simple way.

 We are interested in the exact sequence of cohomology corresponding
to the first of these quadruples.
 Let $L/F$~be a biquadratic field extension with Galois
group~$\Theta$.
 We will consider~$\Theta$ as a two-dimensional vector space over
the field~$\Z/2$.
 Notice that there is a canonical isomorphism between~$\Theta$ and
the dual vector space~$\Theta^*$, given by the only nonzero bilinear
form on~$\Theta$ with the property that $(\theta,\.\theta)=0$
for all $\theta\in\Theta$.
 The three intermediate fields $F\sub K_\theta\sub L$ bijectively
correspond to nonzero elements $a$, $b$, $c$ of the group~$\Theta$
in a natural way.
 The six-term exact sequence of Galois cohomology corresponding
to the exact quadruple~$Q_\Theta$ has the form
 \begin{multline*}
   H^n(L)\op H^n(F)^{\op 3} \overset{(\cor,\.\res^3)}
   \longlongrarrow {\textstyle\bigoplus_{\theta=a,b,c}H^n(K_\theta)}
   \overset{\cor^3} \longrarrow \Theta\ot_{\Z/2} H^n(F)   \\
   \overset{\cup}\lrarrow H^{n+1}(F) \overset{\res}\longrarrow
   H^{n+1}(L) \overset{\cor^3}\longrarrow
   {\textstyle\bigoplus_{\theta=a,b,c}H^{n+1}(K_\theta)},
 \end{multline*}
where the middle arrow $\Theta\ot H^n(F)\rarrow H^{n+1}(F)$
is defined in terms of the embedding $\Theta\simeq
\Hom(\Theta,\.\Z/2)\rarrow \Hom(G_F,\.\Z/2)\simeq H^1(F)$.
 This 6\+term exact sequence is a part of one of the 7\+term
exact sequences of Merkurjev--Tignol and one of the 8\+term
exact sequences of Kahn~\cite{MT,Kahn} for a biquadratic
extension of fields.

\begin{corollary}
 \.Both 8\+term sequences of Kahn~\cite{Kahn} for Galois
cohomology in a biquadratic field extension~$L/F$ are exact,
provided that the first Bockstein homomorphism
$H^n(E,\.\Z/2)\rarrow H^{n+1}(E,\.\Z/2)$ vanishes in the relevant
degree~$n$ for the five fields~$E$ between $F$ and~$L$ (namely,
$F$, $L$, and all three intermediate fields).
\end{corollary}

\begin{proof}
 The middle 6\+term part of the first of Kahn's 8\+term sequences
is written down right above; its exactness follows from Theorem~6
applied to the exact quadruple~$Q_\Theta$ tensored with
$\mu_4^{\ot n}$.
 The middle part of the second sequence can be obtained from
the dual exact quadruple $\Hom_\Z(Q_\Theta,\.\Z)$.
 Exactness at the remaining terms follows from exactness at
the middle terms due to Merkurjev--Tignol and Kahn's results.
\end{proof}

\begin{corollary}
 \.Let $L/F$ be a field extension of degree~$4$ and $M/F$ be its Galois
closure.
 Suppose that the degree $[M:F]$ is a power of\/~$2$.
 Then the kernel~$J$ of the restriction map $H^*(G_F,\.\Z/2)
\rarrow H^*(G_L,\.\Z/2)$ is generated in degrees $1$ and\/~$2$
as an ideal in $H^*(G_F,\.\Z/2)$.
 Moreover, the kernel of the restriction map $H^*(G_F,\.\Z/2)
\rarrow H^*(G_M,\.\Z/2)$ is also generated in degrees $1$ and\/~$2$.
\end{corollary}

\begin{proof}
 There are only three cases: either~$L/F$ is a cyclic extension,
or a biquadratic one, or it has the form~$L'/F$ for a dihedral
extension~$M/F$ with Galois group~$\Gamma_4$.
 The first case follows from Corollary~14, the second from
Corollary~18, and the third from Proposition~17 (both statements).
 Note that in the biquadratic case the ideal~$J$ is actually
generated by a two-dimensional subspace in $H^1(G_F,\.\Z/2)$.
\end{proof}

\begin{corollary}
 \.For any field~$F$ and any two-dimensional subspace\/
$\Theta\sub H^1(G_F,\.\Z/2)$ the kernel of the multiplication
map\/ $\Theta\ot H^*(G_F,\.\Z/2)\rarrow H^{*+1}(G_F,\.\Z/2)$
is a quadratic module over the algebra~$H^*(G_F,\.\Z/2)$.
\end{corollary}

\begin{proof}
 Clearly, any two-dimensional subspace in~$H^1(F)$ corresponds
to a biquadratic field extension~$L/F$, so we can apply
the 6\+term exact sequence of Corollary~18.
 The kernel in question is isomorphic to the quotient module
of $\bigoplus_\theta H^{\ge1}(K_\theta)$ by the image of
the module $H^{\ge1}(L) \op H^{\ge1}(F)^{\op 3}$.
 Since the first module is quadratic (by Theorem~9) and
the second module is generated in degree~$1$ and~$2$
(by Proposition~12), the quotient module is quadratic.
\end{proof}

\begin{conjecture}
 \.For any field~$F$, the ideal generated by any two-dimensional
subspace\/ $\Theta\sub H^1(G_F,\.\Z/2)$ in the cohomology algebra
$H^*(G_F,\.\Z/2)$ is a Koszul module (see~\cite{Pos}) over that
algebra.
\end{conjecture}

\smallskip

\begin{remark}
 It is possible to extend the first statement of Corollary~19
to arbitrary field extensions of degree~$4$ using the following
exact quadruple of permutational representations of
the symmetric group~${\mathbb S}_4$:
 $$
  \Z \lrarrow \Z[X_4]\op\Z \lrarrow \Z[X_6]\op\Z \lrarrow \Z[X_3],
 $$
where $X_4$ is a four-element set on which ${\mathbb S}_4$~acts
in the standard way, $X_6$~is the set of all two-element subsets
of~$X_4$, and $X_3$~is the quotient of~$X_6$ modulo
the involution~$*$ sending a subset to its complement.
 The first map in this sequence is $(\#X_4,{}-2)$, the third map
sends a two-element subset to its class modulo~$*$ and $1$ to
${}-\#X_3$, and the middle arrow is given by the matrix~$U$ with
component maps $U_{11}\:\Z[X_4]\rarrow\Z[X_6]$ sending each
element of~$X_4$ to the sum of three two-element subsets
containing it, $U_{12}=\pr_1$, $\.U_{21}=\#X_6$, and $U_{22}=2$.
 It is easy to construct a homotopy map~$h$ with $dh+hd=2$.
 In fact, one can generalize the 6-term biquadratic exact
sequences to arbitrary field extensions of degree~$4$ using this
exact quadruple!
\end{remark}

\end{document}